\documentclass[conference, 10pt]{IEEEtran}

\usepackage{cite}
\usepackage{url}
\usepackage{color}
\usepackage{placeins}
\usepackage{float}
\usepackage{tabularx,colortbl}
\usepackage{ifthen}

\usepackage{amsmath,amssymb,amsfonts}
\usepackage{algorithmic}
\usepackage{graphicx}
\usepackage{textcomp}
\usepackage{comment}

\usepackage{makecell}
\usepackage{siunitx}
\usepackage{todonotes}

\newcommand{\matel}[1]{\begin{bmatrix} #1 \end{bmatrix}}  
\newcommand{\vt}[1]{\boldsymbol{#1}}  
\newcommand{\mat}[1]{\boldsymbol{\mathsf{#1}}}  
\newcommand{\tmat}[1]{\boldsymbol{\mathsf{#1}}}  
\newcommand{\mrho}[1]{\rho\!\left( #1 \right)}  
\newcommand{\inner}[2]{\langle #1, #2 \rangle}  
\newcommand{\bigo}[1]{\mathcal{O}\!\left( #1 \right)}  
\newcommand{\op}[1]{\mathcal{#1}}  
\newcommand{\inv}[1]{#1^{-1}}  
\newcommand{\tr}[1]{#1^\mathsf{T}}  
\newcommand{\norm}[1]{\left\| #1 \right\|}  

\newcommand{\p}{\vt{r}}  
\newcommand{\unitvec}{\hat{\vt{n}}}  
\newcommand{\de}{\mathrm{d}}  
\newcommand{\eye}{\mat{I}}  

\newcommand{\contr}{\bar{\bar{\chi}}}  
\newcommand{\curr}{\vt{J}}  
\newcommand{\pcurr}{\curr\!_p}  
\newcommand{\contincs}{\xi}
\newcommand{\contincw}{\curr\!_\mathrm{eq}}
\newcommand{\pyr}{p}  
\newcommand{\hatf}{\vt{h}}  
\newcommand{\opS}{\op{S}}  
\newcommand{\opD}{\op{D}^*}  
\newcommand{\opSv}{\op{S}_v}  
\newcommand{\opDv}{\op{D}^*_v}  

\newcommand{\G}{\mat{G}}
\newcommand{\Gss}{\G_{ss}}
\newcommand{\Gww}{\G_{ww}}
\newcommand{\D}{\mat{D}^*}
\newcommand{\Dss}{\D_{ss}}
\newcommand{\Dsw}{\D_{sw}}
\newcommand{\Sl}{\mat{S}}
\newcommand{\Sws}{\Sl_{ws}}
\newcommand{\Sww}{\Sl_{ww}}  

\newcommand{\Z}{\mat{Z}}
\newcommand{\Zss}{\Z_{ss}}
\newcommand{\Zsw}{\Z_{sw}}
\newcommand{\Zws}{\Z_{ws}}
\newcommand{\Zww}{\Z_{ww}}
\newcommand{\Zs}{\Z_\mathrm{self}}
\newcommand{\Zc}{\Z_\mathrm{coupl}}
\newcommand{\inc}{\vt{\alpha}}
\newcommand{\incs}{\inc_s}
\newcommand{\incw}{\inc_w}
\newcommand{\rhs}{\vt{v}}
\newcommand{\rhss}{\rhs_s}
\newcommand{\rhsw}{\rhs_w}

\newcommand{\near}{\mat{N}}
\newcommand{\incww}{\vt{x}}
\newcommand{\rhsww}{\vt{b}}

\newcommand{\Dn}{\D_\mathrm{near}}
\newcommand{\Sn}{\Sl_\mathrm{near}}
\newcommand{\apprDn}{\widetilde{\mat{D}}^*_\mathrm{near}}
\newcommand{\apprSn}{\widetilde{\mat{S}}_\mathrm{near}}
\newcommand{\map}{\mat{\Phi}}
\newcommand{\mapinterp}{\tmat{\Lambda}}
\newcommand{\Toepl}{\tmat{g}}

\makeatletter

\newcounter{author}
\renewcommand{\author}[2][]{
   \stepcounter{author}
   \@namedef{author@\theauthor}{#2}
   \@namedef{authorlabel@\theauthor}{#1}
}

\newcounter{address}
\newcommand{\address}[2][]{
   \stepcounter{address}
   \@namedef{address@\theaddress}{#2}
   \@namedef{addresslabel@\theaddress}{#1}
}

\newcommand{\alsep}{and}

\def\newmaketitle{\par%
  \begingroup%
  \normalfont%
  \def\thefootnote{}
  \def\footnotemark{}
  \let\@makefnmark\relax
  \footnotesize
  \footnotesep 0.7\baselineskip
  \normalsize%
  \twocolumn[\thenewmaketitle\@IEEEaftertitletext]%
  \if@IEEEusingpubid
     \enlargethispage{-\@IEEEpubidpullup}%
  \fi
  \endgroup
  \setcounter{footnote}{0}\let\maketitle\relax\let\@maketitle\relax
  \gdef\@thanks{}%
  \let\thanks\relax}

\def\thenewmaketitle{
  \newpage
  \begin{center}%
    \vskip0.2em{\Huge\@IEEEcompsoconly{\sffamily}\@IEEEcompsocconfonly{\normalfont\normalsize\vskip 2\@IEEEnormalsizeunitybaselineskip
   \bfseries\large}\@title\par}\vskip1.0em\par%
    \vspace{1ex}
    \newcounter{c@author}
    \newcounter{c@tmp}
    \ifthenelse{\value{author}=2}{%
      \newcommand{\liand}{ and }}{%
      \newcommand{\liand}{, and }}
    \ifthenelse{\value{address}<2}{%
      \@nameuse{author@1}%
      \stepcounter{c@author}%
      \whiledo{\value{c@author}<\value{author}}{%
        \setcounter{c@tmp}{\value{author}}%
        \addtocounter{c@tmp}{-\value{c@author}}%
        \ifthenelse{\value{c@tmp}=1}{%
          \renewcommand{\alsep}{\liand}}{\renewcommand{\alsep}{, }}%
        \stepcounter{c@author}\alsep \@nameuse{author@\thec@author}}\\%
    }
    {
      \@nameuse{author@1}${}^{(\ref{\@nameuse{authorlabel@1}})}$%
      \stepcounter{c@author}%
      \whiledo{\value{c@author}<\value{author}}{%
      \setcounter{c@tmp}{\value{author}}%
      \addtocounter{c@tmp}{-\value{c@author}}%
      \ifthenelse{\value{c@tmp}=1}{%
        \renewcommand{\alsep}{\liand}}{\renewcommand{\alsep}{, }}%
      \stepcounter{c@author}\alsep \@nameuse{author@\thec@author}%
        ${}^{(\ref{\@nameuse{authorlabel@\thec@author}})}$%
      }
    }
    \vspace{0.2ex}

    \ifthenelse{\value{address}>0}{%
      \ifthenelse{\value{address}=1}{
        {\@nameuse{address@1}}
      }
      {
        \newcounter{c@address}

        \begin{center}
        \whiledo{\value{c@address}<\value{address}}
        {
          \refstepcounter{c@address}
            ${}^{(\thec@address)}$\,%
              \label{\@nameuse{addresslabel@\thec@address}}%
              \@nameuse{address@\thec@address}\\ %
        }
        \end{center}
      } 
    }
    {
      \relax
    }
  \end{center}
}

\makeatother

\title{On a Fast Solution Strategy for a Surface-Wire Integral Formulation of the Anisotropic Forward Problem in Electroencephalography}
\author[polito]{Carlo Baronio$^{(1)}$}
\author[polito]{Giulio Cosentino$^{(1)}$}
\author[polito]{Paolo Ricci}
\author[imt]{Clément Henry}
\author[polito]{Maxime Y. Monin}
\author[imt]{Adrien Merlini}
\author[polito]{Francesco P. Andriulli}

\address[galfer]{Galileo Ferraris High School, Turin, Italy}
\address[polito]{Politecnico di Torino, Turin, Italy}
\address[imt]{IMT Atlantique, Brest, France}

\begin{document}

\newmaketitle

\begin{abstract}
This work focuses on a quasi-linear-in-complexity strategy for a hybrid surface-wire integral equation solver for the electroencephalography forward problem. The scheme exploits a block diagonally dominant structure of the wire self block---that models the neuronal fibers self interactions---and of the surface self block---modeling interface potentials. This structure leads to two Neumann iteration schemes further accelerated with adaptive integral methods. The resulting algorithm is linear up to logarithmic factors. Numerical results confirm the performance of the method in biomedically relevant scenarios.
\end{abstract}

\section{Introduction}

Several neuro-pathologies require precise functional brain imaging as part of their diagnostic or therapeutic protocols (see~\cite{hybrid} and references therein). Among non-invasive strategies, high resolution electroencephalography (HR-EEG), that images the electric activity of the brain from scalp potentials, is widely used.
In HR-EEG the volume currents are retrieved from the measurements of the electric potentials on the scalp by solving the EEG inverse problem.
Solving this inverse problem requires multiple solutions of the EEG forward problem (FP) in which the surface potential generated by a known current configuration is computed.
Boundary element methods (BEMs) are very popular in the biomedical community to model the FP and a recent hybrid formulation~\cite{hybrid} has introduced the possibility of modeling white matter anisotropies by coupling surface BEM with an integral equation for partially conducting wires.
In this work we present a fast matrix-vector multiplication algorithm for this hybrid formulation which, by exploiting the block diagonal dominance structure (induced by the presence of neuronal fibers in the model) and coupling this matrix structure with adaptive integral methods, obtains a scheme with $O(N \log N)$ complexity in the $N$ degrees of freedom. Theoretical and algorithmic considerations will be complemented by numerical experiments showing the impact of the formulation on medical scenarios.

\section{Background and Notation}

Consider a sequence of nested compartments $\Omega_i$, $i=1,\ldots,C$  modeling the different layers of the head medium, characterized by homogeneous and isotropic conductivities $\sigma_i$. The boundary of each compartment is denoted by  $\Gamma_i$.
Following the strategy in~\cite{hybrid}, the inhomogeneity and anisotropy of the head medium is modeled by populating the white matter with wires of finite anisotropic conductivity contrast $\contr(\p) = \left( \sigma_{i_w} \eye-\bar{\bar{\sigma}}(\p) \right)\inv{\bar{\bar{\sigma}}}(\p)$ with respect to the background conductivity $\sigma_{i_w}$ of the white matter's compartment.
In this setting, the EEG FP consists in finding the electric potential $\phi(\p)$ on the scalp surface $\Gamma_C$ generated by a primary current $\pcurr(\p)$. To do so, the surface $\contincs$ and wire $\contincw$ unknowns (see~\cite{hybrid} for their physical definition) are expanded with discrete basis functions, i.e. $\contincs \approx \sum_{i=1}^{N_s} (\incs)_i \pyr_i$ and $\contincw \approx \sum_{i=1}^{N_w} (\incw)_i \hatf_i$ where $\pyr_i$ and $\hatf_i$ are the 2D and 1D linear Lagrange interpolants, respectively. 
Following a Galerkin approach leads to a linear system of $N = N_s + N_w$ unknowns
\begin{equation}
\label{eqn:hybsys}
    \Z \inc = \matel{\rhss \\ \rhsw} \text{ with } \Z=\matel{-\Gss + \Dss & -\Dsw \\
    -\Sws & \Gww + \Sww}\,,
\end{equation}
and where $\inc = \tr{[\incs\, \incw]}$,
$(\Gss)_{ij} = \tfrac{\sigma_{i_n}+\sigma_{i_n+1}}
                     {2(\sigma_{i_n+1}-\sigma_{i_n})}
               \inner{\pyr_i}{\pyr_j}_\Gamma$,
$(\Dss)_{ij} = \inner{\pyr_i}{\opD \pyr_j}_\Gamma$,
$(\Dsw)_{ij} = \tfrac{1}{\sigma_{i_n}}
               \inner{\pyr_i}{\opDv \contr \hatf_j}_\Gamma$,
$(\Sws)_{ij} = \inner{\hatf_i}{\nabla \opS \pyr_j}_\Omega$,
$(\Gww)_{ij} = \inner{\hatf_i}{\inv{\contr} \hatf_j}_\Omega$,
$(\Sww)_{ij} = \tfrac{1}{\sigma_{i_n}}
               \inner{\hatf_i}{\nabla \opSv \hatf_j}_\Omega$,
$(\rhss)_i = -\tfrac{1}{\sigma_s}
             \inner{\pyr_i}{\opDv \pcurr}_\Gamma$,
and
$(\rhsw)_i = -\tfrac{1}{\sigma_p}
             \inner{\hatf_i}{\nabla \opSv \pcurr}_\Omega$,
with 
$(\op{S}f)(\p) = \int_S G(\p,\p') f(\p')\,\de S'$, $(\op{D}^*f)(\p) = \int_S \unitvec\cdot\nabla G(\p,\p') f(\p')\,\de S'$,
$(\op{S}_v\vt{f})(\p) = \int_V G(\p,\p') \nabla' \cdot \vt{f}(\p')\,\de V'$, and
$(\op{D}^*_v\vt{f})(\p) = \int_V \unitvec\cdot\nabla G(\p,\p') \nabla' \cdot \vt{f}(\p')\,\de V'$.
Above  $\unitvec$ denotes the unit normal vector pointing outwards  $\Gamma_i$ and  $G(\p,\p') = \frac{1}{4\pi\norm{\p-\p'}}$ is the static Green function.
Once \eqref{eqn:hybsys} is solved, $\opS$ and $\opSv$ can be applied to $\inc$ to get the potential $\phi(\p)$ on $\Gamma_C$.

\section{A Fast Solution Strategy}

With respect to a standard integral formulation for isotropic media, corresponding to the left diagonal block in \eqref{eqn:hybsys}, the inclusion of the white matter anisotropy adds a new wire-wire diagonal block and two coupling blocks in the system. First, the new scheme aims at decoupling the surface and wire solution via block diagonal inversion and Neumann series solution of the remainder: after separating diagonal and off-diagonal blocks: $\mat{Z}=\Zs + \Zc$ with $\Zs=\left[\Zss,\mat 0;\mat 0,\Zww\right]$ and $\Zc=\left[\mat 0,\Zsw;\Zws,\mat 0\right]$,
we solve \eqref{eqn:hybsys} as $\big(\eye + \inv{\Zs}\Zc\big) \inc = \inv{\Zs}\rhs$ via  a Neumann series approach enabled by the block diagonal dominance, in cases of practical relevance, of the original matrix (i.e. for the spectral radius $\rho_Z = \mrho{\inv{\Zs}\Zc}<1$). Thus we have
$\incs^{(k+1)} = \inv{\Zss}(\rhss - \Zsw\incw^{(k)})$ and $\incw^{(k+1)} = \inv{\Zww}(\rhsw - \Zws\incs^{(k)})$
whose complexity reduces to the one of the two inversions and of the multiplication of the coupling terms.
 The  multiplication of the coupling terms can be done efficiently if a fast matrix vector product algorithm is available. We have opted for an adaptive integral method (AIM)~\cite{Bleszynski}.
  In other words, all kernel interactions in $\Dss$, $\Sww$, $\Sws$ and $\Dsw$
between all Gaussian quadrature points are interpolated on the same Cartesian grid with a number of nodes proportional to the number of unknowns $N$ and handled via FFT in $\bigo{N \log N}$ complexity. As is standard in AIM~\cite{Bleszynski}, a near field precorrection is required for all kernels:
a generic $\D$ and $\Sl$ (for surface, wire, and off diagonal couplings) is written as $\D = \Dn - \apprDn + \map_p \tr{\mapinterp} \Toepl_D \mapinterp \map_f$
and $\Sl = \Sn - \apprSn + \map_p \tr{\mapinterp} \Toepl_S \mapinterp \map_f$,
where $\Dn$ and $\Sn$ are the uncompressed near fields, $\apprDn$ and $\apprSn$ are the FFT precorrections, $\mapinterp$ is unique for every product and interpolates the quadrature points, and $\map_p$ and $\map_f$ map quadrature points to basis functions; all these matrices are sparse. The FFT is applied to the Toeplitz matrices $\Toepl_S$ and $\Toepl_D$ that, because of the translation invariance of all Green functions involved, require $\bigo{N}$ memory storage.
Since the double layer kernel is $\unitvec \cdot \nabla G(\p,\p') = \unitvec \cdot \frac{\p'-\p}{4\pi\norm{\p'-\p}^3}$, the product of $\Toepl_D$ with a vector is split into three scalar components.

Since the $\Zss$ block corresponds to the classical homogeneous multilayer BEM formulation, once a fast matrix vector product algorithm is available, it can be inverted iteratively with standard techniques (see \cite{hybrid} and references therein). Regarding $\Zww$, the near field kernel interactions are extracted with an octree and the resulting sparse matrix $\near$ is used as a preconditioner of the linear system $\Zww\incww = \rhsww$.
The near field dominance of $\Zww$---due to the electric current flowing along the fibers, i.e. $\rho_w = \mrho{\inv{\near}(\Zww-\near)} < 1$---enables a second usage of a Neumann series from which
$\incww^{(k+1)} = \incww^{(k)} + \inv{\near}\big(\rhsww - \Zww\incww^{(k)}\big)$.
A sparse solver is used to invert $\near$ in $\bigo{N}$ time complexity and the multiplication of $\Zww$ is done in $\bigo{N \log N}$ with the AIM.

\begin{table}
\caption{Performance Comparison}
\centering
\begin{tabular}{c|c|c|c}
\hline
    & \bf Setup time & \bf Storage & \bf Time per RHS \\
    \hline
    \makecell{Standard Iterative Solution} & \SI{24657}{\second} & \SI{31.2}{GB} & \SI{89.15}{\second}  \\ 
    \makecell{This work} & \SI{2498}{\second} & \SI{0.4}{GB} & \SI{26.42}{\second} \\ 
    \hline
\end{tabular}
\vspace{-0.2cm}
\label{tab}
\end{table}

\begin{figure}
\centerline{\includegraphics[scale=0.45]{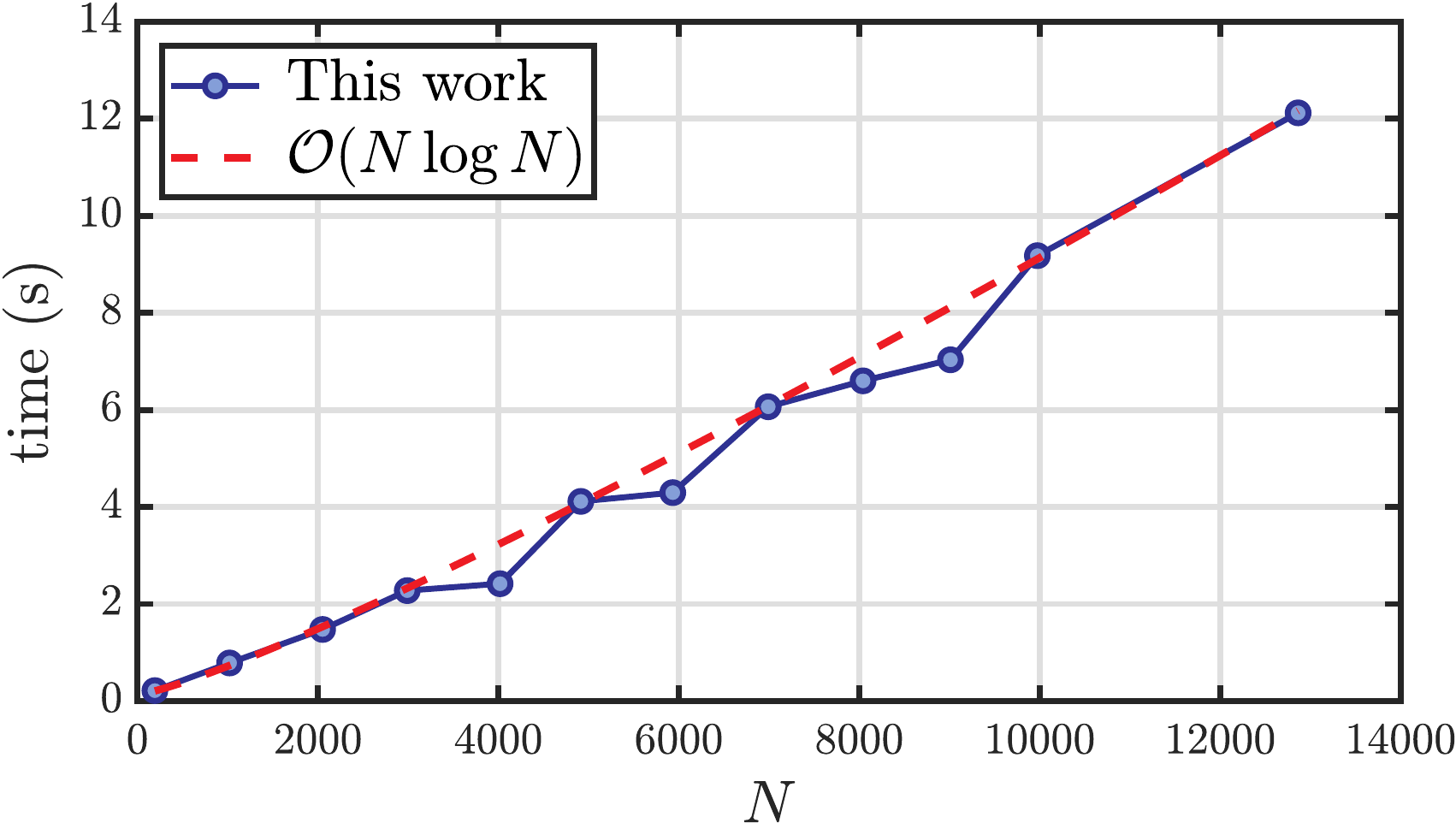}}
\caption{Solution time vs number of unknowns.}
\vspace{-0.5cm}
\label{time-sol}
\end{figure}

\section{Numerical Results}
\begin{figure}
\centerline{\includegraphics[scale=0.2]{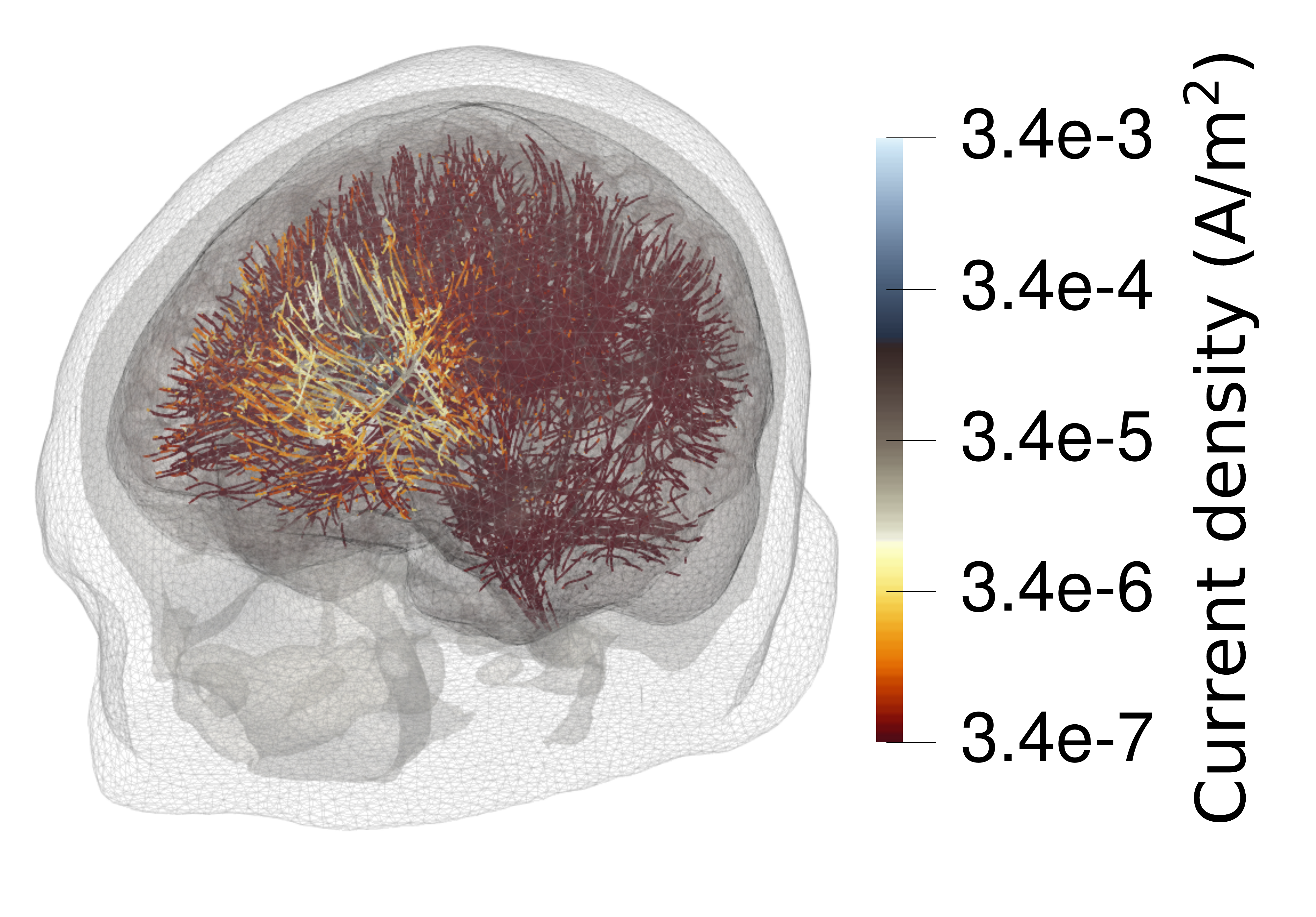}} 
\vspace{-0.3cm}
\caption{Brain fibers current density on the MRI head model.}
\label{brain}
\vspace{-0.4cm}
\end{figure}

The favorable complexity scaling of the proposed scheme has been verified on a set of canonical geometries composed of spherical surfaces and
orthogonal brain fibers. The total timings are reported in Fig.~\ref{time-sol} and clearly confirm that the scheme we propose is, up to logarithmic factors, linear in complexity.
The relevance of our fast solution strategy for real case scenarios has been studied on a realistic head model obtained with magnetic resonance imaging (MRI) data that includes white matter neuronal fibers with a tangential anisotropic conductivity of $\SI{1.3}{\siemens\per\meter}$ and four layers (gray matter, cerebrospinal fluid, skull, scalp) with conductivities $\SI{0.13}{\siemens\per\meter}$, $\SI{1.79}{\siemens\per\meter}$, $\SI{0.01}{\siemens\per\meter}$, and $\SI{0.43}{\siemens\per\meter}$ respectively. The obtained current on the neuronal fibers is shown in Fig.~\ref{brain}. 
For this problem the radius of the fibers is chosen to match a total volume of $\SI{450}{\cubic\milli\meter}$. The total number of unknowns is $63\,922$. The two spectral radii are $\rho_Z = 0.439$ and $\rho_w = 0.799$, both less than one, thus allowing the Neumann strategy.
For this experiment we have compared in Table~\ref{tab} the method proposed in this work
with the uncompressed solution. In both cases the tolerance iterative schemes has been set to $10^{-3}$ and the results show the advantage of the new scheme. 

\section*{Acknowledgment and Contributions}\label{sec:ack}
The work of this paper and the associated Early-Research Program for talented high school students has received funding from the European Research Council (ERC) under the European Union’s Horizon 2020 research and innovation programme (grant agreement No 724846, project 321) and from the EU H2020 research and innovation programme under the Marie Skłodowska-Curie grant agreement n° 955476 (project COMPETE). The first two authors are listed in alphabetic order and have contributed equally to this work.



\begin{thebibliography}{1}

\bibitem{hybrid} M. Y. Monin, L. Rahmouni, A. Merlini and F. P. Andriulli, ``A hybrid volume-surface-wire integral equation for the anisotropic forward problem in electroencephalography,'' IEEE J. Electromagn. RF Microw. Med. Biol., vol. 4, no. 4, pp. 286--293, Dec. 2020.
\bibitem{Bleszynski} E. Bleszynski, M. Bleszynski and T. Jaroszewicz, ``AIM: Adaptive integral method for solving large-scale electromagnetic scattering and radiation problems,'' in Radio Science, vol. 31, no. 5, pp. 1225--1251, Sept.-Oct. 1996.

\end{thebibliography}
%

\end{document}